\documentclass[aps,pre,reprint,groupedaddress]{revtex4-2}

\usepackage{amsmath}
\usepackage{amssymb}
\usepackage{graphicx}
\usepackage{xspace}
\usepackage{tabularx}
\usepackage{color}
\usepackage{cases}
\usepackage{ragged2e}
\usepackage{listings}
\usepackage{booktabs}
\usepackage{colortbl}

\newcommand{\be}{\begin{equation}}
\newcommand{\ee}{\end{equation}}
\newcommand{\bea}{\begin{subequations}\begin{eqnarray}}
\newcommand{\eea}{\end{eqnarray}\end{subequations}}
\newcommand{\eq}[1]{Equation~(\ref{eq:#1})}
\newcommand{\eqs}[2]{Equations~(\ref{eq:#1}) and (\ref{eq:#2})}
\newcommand{\eqp}[1]{(Equation~\ref{eq:#1})}
\newcommand{\fig}[1]{Figure~\ref{fig:#1}}
\newcommand{\tab}[1]{Table~\ref{tab:#1}}
\newcommand{\s}[1]{Section~\ref{s:#1}}
\renewcommand{\ss}[1]{Section~\ref{ss:#1}}

\renewcommand{\ap}[1]{Appendix~\ref{a:#1}}

\newcommand{\msun}{{\rm M}_{\odot}}

\newcommand{\eg}{e.g.,\xspace}
\newcommand{\ie}{i.e.,\xspace}
\newcommand{\cf}{c.f.,\xspace}
\newcommand{\R}{\mathbb{R}}
\newcommand{\N}{\mathbb{N}}
\newcommand{\C}{\mathbb{C}}
\renewcommand{\d}{{\rm d}}

\renewcommand{\L}{\mathcal{L}}
\newcommand{\erf}{{\rm erf}}
\newcommand{\erfi}{{\rm erfi}}
\def\code#1{\texttt{#1}}
\newcommand{\Rcode}{\textsc{R}\xspace}
\newcommand{\Pythoncode}{\textsc{Python}\xspace}

\begin{document}


\title{From Cavitation to Astrophysics: Explicit Solution of the Spherical Collapse Equation}


\author{Danail Obreschkow\vspace{2mm}}
\affiliation{International Centre for Radio Astronomy Research (ICRAR) and International Space Centre (ISC), M468, University of Western Australia, 35 Stirling Hwy, Perth, WA 6009, Australia\vspace{2mm}}


\date{Accepted by Physical Review E on 16 May 2024}

\begin{abstract}
Differential equations of the form $\ddot R=-kR^\gamma$, with a positive constant $k$ and real parameter $\gamma$, are fundamental in describing phenomena such as the spherical gravitational collapse ($\gamma=-2$), the implosion of cavitation bubbles ($\gamma=-4$) and the orbital decay in binary black holes ($\gamma=-7$). While explicit elemental solutions exist for select integer values of $\gamma$, more comprehensive solutions encompassing larger subsets of $\gamma$ have been independently developed in hydrostatics (see Lane-Emden equation) and hydrodynamics (see Rayleigh-Plesset equation). I here present a universal explicit solution for all real $\gamma$, invoking the beta distribution. Although standard numerical ODE solvers can readily evaluate more general second order differential equations, this explicit solution reveals a hidden connection between collapse motions and probability theory that enables further analytical manipulations, it conceptually unifies distinct fields, and it offers insights into symmetry properties, thereby enhancing our understanding of these pervasive differential equations.
\end{abstract}


\maketitle



\section{Introduction}\label{s:introduction}

It may be surprising that the 21$^{\rm st}$ century still holds critically important differential equations, which admit compact, but barely known explicit solutions. In this paper, I discuss such a family of equations that govern all physical systems described by a one-dimensional time-dependent coordinate $R(T)$, subject to a restoring force that varies as a power law of that coordinate. Formally, this is expressed by the non-linear ordinary differential equation (ODE)
\be\label{eq:main}
	\ddot R = -kR^\gamma
\ee
with fixed $k\in\R_+$ and $\gamma\in\R$. Dots denote derivatives with respect to time $T$. Using a capital $T$, normally reserved for temperature, will help distinguish this time from a dimensionless time $t$ (\ss{standardform}). We may call \eq{main} the `spherical collapse equation' by virtue of its most common applications.

I will restrict the discussion to the case, where $R$ exhibits a finite maximum $R_0$---in many physical applications this is equivalent to stating that the system is `bound'. \eq{main} being an autonomous ODE, we can always choose the origin of time to coincide with $R=R_0$, in which case the initial conditions read
\be\label{eq:main.ics}
	R(0)=R_0~~~\text{and}~~~\dot R(0)=0.
\ee

\eqs{main}{main.ics} are manifestly time-symmetric, \ie $R(T)=R(-T)$ where a solution exists. Often most relevant is the time interval $T\in[0,T_c]$, from the maximum radius to the (first) collapse point, $R(T_c)=0$. I will focus on this interval for most of this paper.

Table~\ref{tab:applications} lists textbook examples of physical systems governed by \eq{main}. Many of them are encountered at undergraduate levels and used to illustrate cases, whose solution $R(T)$ must be obtained via numerical ODE solvers, such as Runge--Kutta solvers. The existence, let alone the form, of an explicit solution is rarely noted, except, of course, in trivial cases like the linear free-fall ($\gamma=0$) and harmonic oscillator ($\gamma=1$).

An iconic example of \eq{main} in astrophysics is the gravitational collapse of a uniform collisionless sphere ($\gamma=-2$), the so-called top-hat spherical collapse model \citep{Gunn1972}. A simple parametric solution, $\{T(\theta),R(\theta)\}$ \citep{Lin1965}, has become the default solution in cosmology. Its transcendental nature implies that algebraic solutions of \eq{main} do not generally exist, yet the search for non-algebraic \cite{Slepian2023} and approximate \citep{Foong2008} solutions continues.

Barely known to astrophysicists, an analogous collapse equation (but with $\gamma=-4$) has long been studied in hydrodynamics \cite{Stokes1847}. This equation, named after Lord Rayleigh \cite{Rayleigh1917}, describes the collapse of an empty spherical cavity in an ideal incompressible liquid. Recent approximations \cite{Obreschkow2012a} have led to the development of an infinite series that rapidly converges toward the exact solution \cite{Amore2013}. Shortly after, others presented the first closed-form solution for cavities in three \citep{Kudryashov2014} and $N\geq3$ \citep{Kudryashov2015} dimensions, in terms of hypergeometric functions.

\eq{main} is also encountered in static systems. For example, it is equivalent to the one-dimensional Lane-Emden equation \citep{Lane1870,Emden1907}, describing the density profile of a self-gravitating polytropic gas in a thin tube in hydrostatic equilibrium. A classic solution for $\gamma\!>\!-1$ has been found in the context of galactic discs \cite{Harrison1972}.

Building on the special solutions found independently in different fields, this paper presents a compact general explicit solution of \eq{main}. While overall straightforward, the path to this solution involves some subtleties that enable the generalisation to all real $\gamma$. This allows us to discuss various seemingly unrelated basic problems in physics and astrophysics in unison. \s{background} overviews reformulations, useful for the derivation of the general solution in \s{solution} and its discussion in \s{discussion}. \s{conclusion} concludes with a brief synthesis and outlook.

\begin{table*}
\newcolumntype{c}[1]{>{\centering\arraybackslash}p{#1}}
\newcolumntype{l}[1]{>{\RaggedRight\arraybackslash}p{#1}}
\newcommand{\q}[1]{\textbf{#1}\hfill\break}
\newcommand{\img}[1]{\raisebox{-1.09cm}{\includegraphics[width=2cm]{schematics_#1}}}
\definecolor{lightgray}{gray}{0.8}
\newcommand{\midrulecol}{\arrayrulecolor{lightgray} \\[-4.3ex]\midrule \arrayrulecolor{black}\\[-4ex]}
\setlength{\tabcolsep}{4.5pt}
\renewcommand{\arraystretch}{1.5}
\begin{tabular}{l{5.5cm}p{2cm}l{2.3cm}l{4.3cm}c{0.6cm}c{1.3cm}}
\toprule
\bf{Physical system} & \bf{Schematic} & \bf{Meaning of} $R$ & \bf{Time scale} $T_0$ & $\gamma$ & $\tau$ \\ [1ex] \toprule
\q{Cavitation bubble} Spherical cavity--empty or with low constant inner pressure--imploding under the pressure of an incompressible inviscid liquid without surface tension & \vspace{-2.5mm}\img{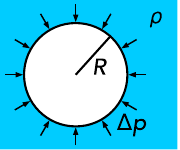} & Bubble radius & $R_0\sqrt{\rho/\Delta p}$, where $\rho$ is the density of the liquid and $\Delta p=p_\infty-p_v$, with far-field pressure $p_\infty$ and constant bubble pressure $p_v\ll p_\infty$ & $-4$ & $\approx0.91468$ \\ \midrulecol
Generalisation to a spherical bubble in $N\geq3$ dimensions \citep[\eg][]{Klotz2013} & & & $(N-2)^{-1/2}R_0\sqrt{\rho/\Delta p}$ & \mbox{$\!\!\!\!\!-N-1$} & ~(Eq.~\ref{eq:collapsetime}) \\ [3.5ex] \midrule
\q{Gravitational collapse} Collapse of a uniform sphere, also known as spherical top-hat, subjected only to gravitational forces & \img{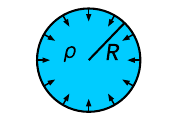} & Radius of sphere  & $\sqrt{3/(4\pi G\rho_0)}$, where $\rho_0$ is the initial density at the onset of the collapse and $G$ is the gravitational constant& $-2$ & $\pi/\sqrt{8}\approx1.11072$ \\ \midrule
\q{Two-body collision} Two point masses in gravitational free-fall toward each other, on a straight line without angular momentum & \img{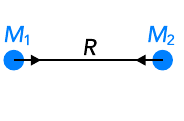} & Distance between the centres of the two masses & $R_0^{3/2}/\sqrt{G(M_1+M_2)}$, where $M_i$ denote the two masses and $G$ is the gravitational constant & $-2$ & $\pi/\sqrt{8}\approx1.11072$ \\ \midrule
\q{Free-fall in spherical potential} Radial free-fall in a spherically symmetric gravitational power law potential, $\phi\sim R^\alpha$, with index $\alpha\ne0$ & \img{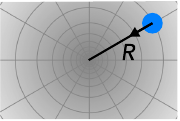} & Distance from centre of potential & $\sqrt{(L/R_0)^\alpha/\alpha}\cdot R_0/V$, where $L$ and $V$ are the length and velocity scales of the potential $\phi(R)={\rm sgn}(\alpha)\cdot V^2\cdot(R/L)^\alpha$ & \mbox{$\alpha-1$} & ~(Eq.~\ref{eq:collapsetime}) \\ \midrulecol
Same for a logarithmic potential, \eg the potential of a singular isothermal sphere\!\!\!\!\!\!\!\!\!\! & & & $R_0/V$, where $V$ is the velocity scale in $\phi(R)=V^2\ln(R/L)$ & $-1$ & $\sqrt{\pi/2}\approx1.25331$ \\ [3ex]  \midrule
\q{Free-fall in uniform field} Straight fall in a constant acceleration field, such as a vertical drag-free drop of an object on Earth & \img{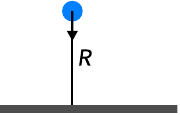} & Height from the ground & $\sqrt{R_0/g}$, where $g$ is the norm of the constant acceleration, \eg $g=9.81\rm m\,s^{-2}$ on Earth & $0$ & $\sqrt{2}\approx1.41421$ \\ \midrule
\q{Acceleration by a dipole charge} Straight fall of a free charge toward a fixed dipole of charges, separated by a distance $\varepsilon\ll R$ along the trajectory& \img{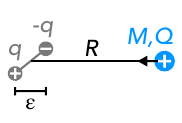} & Distance between free charge and dipole centre & $\sqrt{MR_0^4/(4\varepsilon k_e Qq)}$, where $M$ is the mass of the free charge $Q\!>\!0$, $q\!>\!0$ the dipole charge and $k_e$ the Coulomb constant & $-3$ & $1$ \\ \midrule
\q{Relativistic orbital decay} Orbital decay by gravitational radiation of two masses at non-relativistic velocities, initially on circular orbits & \img{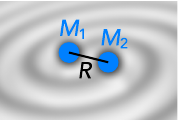} & Distance between the centres of the two masses & $\frac{5c^5R_0^4}{64\sqrt{3}G^3M_1M_2(M_1+M_2)}$, with $M_i$ the two masses, $G$ the gravitational constant and $c$ the speed of light & $-7$ & $\approx0.74683$ \\ \midrule 
\q{Harmonic oscillator} Motion of a mass attached to an ideal spring, initially at rest in an extended or compressed state & \img{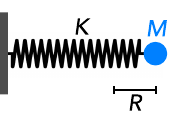} & Distance from the equilibrium position & $\sqrt{M/K}$, where $M$ is the mass and $K$ the spring constant & $+1$ & $\pi/2\approx1.57080$ \\ \midrule
\q{Polytrope} Density profile of a self-gravitating polytropic gas in one dimension in hydrostatic equilibrium & \img{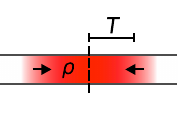} & $\rho^{1/n}$, where $\rho$ is the gas density and $n$ the index of $p=K\rho^{1+1/n}$ & $\sqrt{(1+n)KR_0^{1-n}/(4\pi G)}$ is a length scale set by the parameters of the equation of state $p=K\rho^{1+1/n}$ & $n$ & (Eq.~\ref{eq:collapsetime}) \\
\bottomrule
\end{tabular}
\caption{Textbook examples of systems governed by \eqs{main}{main.ics}. Moving objects described by $R(T)$ are shown in blue, whereas static objects described by this equation are shown in red. $T_0$ is the characteristic physical scale (normally a characteristic time), such that \eq{main} reduces to the dimensionless \eq{main.standard} upon normalisation via \eq{nondimensional}. Where $T_0$ depends on $R_0$, the latter is the maximum, initial value of $R$. In the case of the polytrope (bottom row), $R$ and $T$ have different meanings: $R:=\rho^{1/n}$ is a measure of the gas density and $T$ is the distance from the centre of mass. The last column shows the dimensionless collapse time, $\tau$, such that the physical collapse time is $T_c=\tau T_0$, where $R(T_c)=0$. For most examples, analogous cases with other forces, such as electric forces instead of gravitational ones, can be found.}
\label{tab:applications}
\end{table*}


\section{Setting the stage}\label{s:background}

Before delving into the solutions of \eq{main}, it is worth recalling basic, but pivotal remarks on equivalent reformulations. These will help the derivation (\s{solution}) and discussion (\s{discussion}) of the solutions.

\subsection{Dimensionless form}\label{ss:standardform}

The first remark is that it is convenient to work in dimensionless position and time coordinates, defined as
\begin{equation}\label{eq:nondimensional}
	r := \frac{R}{R_0}~~~{\rm and}~~~t:=\frac{T}{T_0},
\end{equation}
where $R_0=R(0)$ is the maximum, initial value of $R$ (\cf Equation~\ref{eq:main.ics}) and
\be\label{eq:t0}
	T_0:=\sqrt{k^{-1}R_0^{1-\gamma}}
\ee
is the natural time constant. With this normalisation, \eq{main} becomes
\be\label{eq:main.standard}
	\ddot r = -r^\gamma,
\ee
and the initial conditions \eqp{main.ics} become
\be\label{eq:main.ics.standard}
	r(0)=1~~~{\rm and}~~~\dot r(0)=0.
\ee
It is understood that dots above $r(t)$ now denote derivatives with respect to $t$, rather than $T$. To help the following derivations, we introduce the variable
\be\label{eq:nondimensionaltime}
	\tau:=\frac{T_c}{T_0}
\ee
as a shorthand for the dimensionless collapse time.

The constant $k$ has conveniently disappeared in the dimensionless form of \eq{main.standard}, reducing the family of ODEs to a single control parameter $\gamma$. The natural disappearance of $k$ brings to light the inherently scale-free nature of \eq{main} implied by its power law structure. This is the deeper reason for \eq{main} to be applicable from microscopic to astrophysical scales. 

As a textbook example, let us consider the gravitational collapse of a uniform pressure-free sphere, initially at rest with radius $R(0)=R_0$. The Newtonian equation of motion reads
\be\label{eq:tophat}
	\frac{\d^2R}{\d T^2} = -\frac{GM}{R^2},
\ee
which is \eq{main} with $\gamma=-2$ and $k=GM$. Upon expressing $R$ and $T$ in their natural units $R_0$ and $T_0=\sqrt{R_0^3/(GM)}$, \eq{tophat} becomes indeed $\ddot r=-r^{-2}$.


\subsection{Equivalent formulations}

The second remark is that \eq{main.standard}, and thus \eq{main}, can be rewritten in other differential forms.

\subsubsection{Integral of motion}

Most importantly, \eq{main.standard} always exhibits an \emph{integral of motion},
\begin{subnumcases}
	{0=\label{eq:integral}}
	\tfrac{1+\gamma}{2}\dot r^2+r^{1+\gamma}-1 & if $\gamma\ne-1$ \label{eq:integral.a}\\
	\dot r^2+2\ln r & if $\gamma=-1$ \label{eq:integral.b}.
\end{subnumcases}   
If differentiated with respect to $t$ and solved for $\ddot r$, \eq{integral} becomes \eq{main.standard}, except if $\dot r=0$. Since Equations~(\ref{eq:integral}) are first-order differential equations, only one boundary condition is required, such as $r(0)=1$.

\subsubsection{Hybrid differential equation}

If $\gamma\ne-1$, a hybrid ODE, mixing first and second derivatives, can be obtained by rewriting $r^{1+\gamma}$ in \eq{integral.a} as $rr^\gamma=-r\ddot r$ \eqp{main.standard}, hence
\be\label{eq:hybrid}
	1-\tfrac{1+\gamma}{2}\dot r^2+r\ddot r=0.
\ee
This is, for example, the standard form of the equation describing the collapse of a spherical cavitation bubble without viscosity and surface tension \citep{Rayleigh1917}, normally written in dimensional form,
\be\label{eq:rayleigh}
	\frac{3}{2}\left(\frac{\d R}{\d T}\right)^2+R\frac{\d^2R}{\d T^2}+\frac{\Delta p}{\rho}=0,
\ee
where $\Delta p$ is the driving pressure and $\rho$ the constant liquid density (see first row in \tab{applications}). Expressed in natural units, $R_0=R(0)$ and $T_0=R_0\sqrt{\rho/\Delta p}$, \eq{rayleigh} becomes \eq{hybrid} with $\gamma=-4$. Thus \eq{rayleigh} is equivalent to \eq{main} with $\gamma=-4$ and $k=R_0^3\Delta p/\rho$.

\begin{figure}[b]
    \centering
    \includegraphics[width=\columnwidth]{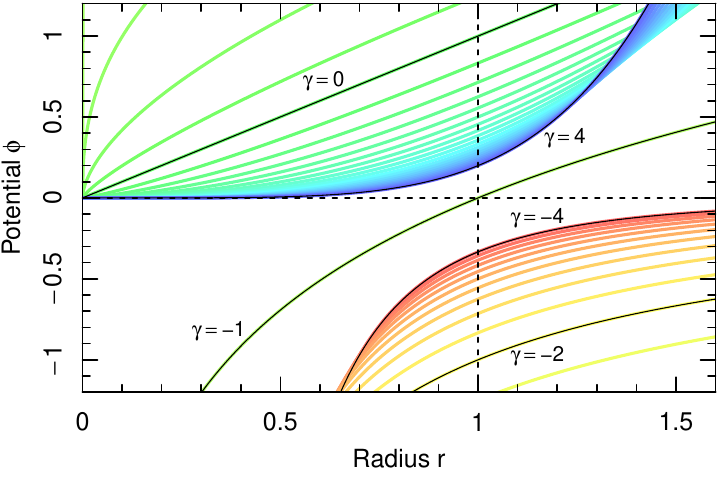}\vspace{-2mm}
    \caption{Potential $\phi(r)$, given in \eq{potential}, for different values of $\gamma$, equally spaced by $\Delta\gamma=0.2$. There are qualitatively distinct families at $\gamma<-1$ (red--yellow) and $\gamma>-1$ (blue--green), as discussed in \ss{collapsepoint}. They are separated by an isolated critical curve for $\gamma=-1$. Thin black lines highlight the coloured curves with a printed $\gamma$-value.}
    \label{fig:potential}
\end{figure}

\subsubsection{Lagrangian}

Finally, with an eye on applications in physics, it is worth noting that \eq{main.standard} derives from a time-independent Lagrangian
\be\label{eq:lagrangian}
	\L(r,\dot r) = \frac{\dot r^2}{2}-\phi(r)
\ee
with the potential, defined up to an additive constant,
\be\label{eq:potential}
	\phi(r)=\left\{\begin{array}{ll}
	r^{1+\gamma}/(1+\gamma) & \rm{if}~\gamma\ne-1 \\
	\ln r & \rm{if}~\gamma=-1
	\end{array}\right..
\ee
The Euler-Lagrange equation, $\tfrac{\d}{\d t}(\partial\L/\partial\dot r)=\partial\L/\partial r$, then generates \eq{main.standard}, for all $\gamma$.

Selected potentials $\phi(r)$ are shown in \fig{potential}. The existence of two regimes, separated by the critical value $\gamma=-1$, is the fundamental reason for the qualitatively different behaviour of $r(t)$ in these regimes. These regimes will be discussed in more detail in \ss{collapsepoint}.


\section{General solution}\label{s:solution}

In this section, I present compact general solutions of \eq{main} with initial conditions of \eq{main.ics}, following an approach similar to \cite{Kudryashov2014}, but for general real parameters $\gamma$. For convenience, all derivations and solutions are presented for the dimensionless form in \eq{main.standard} with initial conditions stated in \eq{main.ics.standard}. Any solution $r(t)$ can be transformed to the corresponding dimensional solution $R(T)$ via the linear transformations of \eq{nondimensional}.

\subsection{Inverse solution $t(r)$}

\subsubsection{Cases $\gamma\ne-1$}

This case can be solved, starting from the integral of motion in \eq{integral.a}, by reduction to quadrature, a standard method expressing the solution of ODEs in terms of integrals. For general first order ODEs, $\dot r(t)=f(r(t), t)$, the integrals found in this way must be evaluated numerically. Interestingly, however, the limited form of \eq{main.standard} allows us, through adequate substitutions, to express the emerging integrals through the \emph{incomplete beta function} $B(x;\alpha,\beta)$ (derivation in \ap{derivation}). This special function benefits from an extensive mathematical literature available for further analytical manipulations of the solution. Moreover, by introducing the shorthands
\be\label{eq:defeta}
	\eta=\frac{1}{|1+\gamma|},
\ee
and
\be\label{eq:defalpha}
	\alpha = \frac{1}{4}+\frac{3-\gamma}{4|1+\gamma|},
\ee
it turns out possible to compact the explicit solution $r(t)$ to a universal form, valid for all $\gamma\ne1$, hence simplifying and generalising state-of-the-art hypergeometric solutions for $\gamma=-4$ \cite{Kudryashov2014} and lower negative integer $\gamma$ \cite{Kudryashov2015}.

Explicitly, \eqs{main.standard}{main.ics.standard} solve to
\be\label{eq:timesolution1}
	t(r) = \tau-\sqrt{\frac{\eta}{2}}~B\left(r^{|1+\gamma|};\alpha,\tfrac{1}{2}\right),
\ee
following the derivation provided in \ap{derivation}. The collapse time $\tau$ is obtained by solving \eq{timesolution1} for $\tau$ at $(r,t)=(1,0)$,
\be\label{eq:collapsetime}
	\tau(\gamma)=\sqrt{\frac{\eta}{2}}~B\big(\alpha,\tfrac{1}{2}\big).
\ee
Here $B(\alpha,\beta)$ is the (complete) \emph{beta function}, \ie the incomplete $B(x;\alpha,\beta)$ evaluated at $x=1$. Selected explicit values of $\tau$ are listed in \tab{tau}.

Substituting \eq{collapsetime} back into \eq{timesolution1}, the latter can be transformed to the special solution for $\gamma>-1$ found in modelling polytropic gas densities in galactic disks \cite{Harrison1972}.
%
%

\subsubsection{Special case $\gamma=-1$}

In this case, \eq{integral.b} integrates to
\be\label{eq:timesolutionspecial}
	t(r) = \tau(-1)\,\erf\Big(\sqrt{-\ln r}\,\Big)
\ee
where
\be\label{eq:collapsetimespecial}
	\tau(-1) = \sqrt{\pi/2}
\ee
is the collapse time, reached for $r\rightarrow0_+$. This collapse time is equal to the limit of \eq{collapsetime} for $\gamma\rightarrow-1$. Likewise, \eq{timesolutionspecial} is the limit of \eq{timesolution1} for $\gamma\rightarrow-1$, showing that these singularities of \eqs{timesolution1}{collapsetime} at $\gamma=-1$ are removable.

\subsection{Explicit solution $r(t)$}

To invert the solutions $t(r)$ it is convenient to normalise \eq{timesolution1} by $\tau$, which yields
\be\label{eq:timesolution2}
	\frac{t}{\tau} = 1-I\left(r^{|1+\gamma|};\alpha,\tfrac{1}{2}\right),
\ee
where $I(x;\alpha,\beta)\equiv B(x;\alpha,\beta)/B(\alpha,\beta)$ is called the \emph{regularised incomplete beta function} (a simple, but useful step pointed out by Roberto Iacono, priv.~com.). \eqs{timesolution2}{timesolutionspecial} can be readily inverted to a compact general explicit solution in terms of well-known special functions,
\begin{subnumcases}
	{r(t) = \label{eq:radiussolution}}
	Q\left(1-\frac{|t|}{\tau(\gamma)};\alpha,\frac{1}{2}\right)^\eta & if $\gamma\!\ne\!-1$ \label{eq:radiussolution.a}\\
	\exp\left(-\erfi^2\Big(\sqrt{2/\pi}~t\Big)\right) & if $\gamma\!=\!-1$~,~~~~ \label{eq:radiussolution.b}
\end{subnumcases}
where $\erfi(x)$ is the \emph{inverse error function} and $Q(x;\alpha,\beta)$ is the \emph{inverse regularised incomplete beta function}. \eq{radiussolution.b} is the limit of \eq{radiussolution.a} for $\gamma\rightarrow-1$, again showing that this singularity of \eq{radiussolution.a} is of a removable type.

The absolute value $|t|$ in \eq{radiussolution.a} does not follow from \eq{timesolution2}, only valid for $t\in[0,\tau]$. It is an ad-hoc extension of $r(t)$ to the domain $t\in[-\tau,\tau]$, exploiting the time-reversal symmetry of \eqs{main.standard}{main.ics.standard}. \eq{radiussolution.b} is already time-symmetric and hence valid on $t\in[-\tau,\tau]$.

The significance of \eq{radiussolution} relies in the fact that the regularised incomplete beta function $I(x;\alpha,\beta)$ is the cumulative density function of the \emph{beta distribution}, a common distribution function in probability theory, related to, but not to be confused with, the beta function. Hence, the function $Q(x;\alpha,\beta)$ is the quantile function of the beta distribution. This connection to the beta distribution, which has been intensively studied in probability theory, makes \eq{radiussolution} a practical solution for further analytical manipulations. Moreover, the quantile function $Q(x;\alpha,\beta)$ is readily accessible to all modern programming languages, making this solution easy to implement. Two examples in \Pythoncode and \Rcode code are given in \ap{implementation}.

\begin{figure}
    \centering
    \includegraphics[width=\columnwidth]{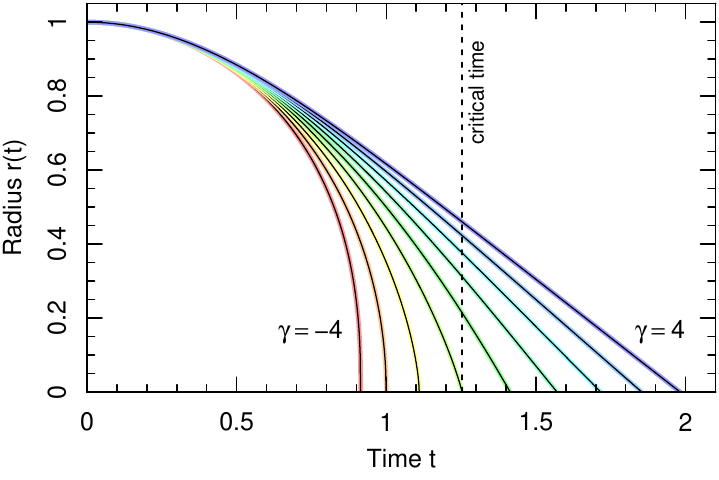}\vspace{-2mm}
    \caption{Exact solutions of the normalised collapse \eq{main.standard} with initial conditions of \eq{main.ics.standard}, evaluated on $t\in[0,\tau(\gamma)]$ using the explicit \eq{radiussolution}. Coloured curves correspond to different values of $\gamma$, separated by $\Delta\gamma=1$. The dashed vertical line marks the collapse time $\sqrt{\pi/2}$, corresponding to the critical exponent $\gamma=-1$ that separates the divergent behaviour of $\dot r(\tau)$ from the convergent one.}
    \label{fig:rt}
\end{figure}

\fig{rt} shows $r(t)$ on the interval $t\in[0,\tau(\gamma)]$ for all integer $\gamma$ from $-4$ to $+4$. Key properties of these solutions will be discussed in \s{discussion}.

\begin{table}
\centering
\begin{tabularx}{\columnwidth}{@{\extracolsep{\fill}}ccccc}
\toprule
$\gamma$ & $\tau$ (exact) & $\tau$ (num) & $\dot r(\tau)$ (exact) & $\dot r(\tau)$ (num) \\ [0.5ex]
\midrule \\ [-2ex]
$-\infty$ & $0$ & 0 & $-\infty$ & $-\infty$ \\
$-100$ & --- & $0.22019512$ & $-\infty$ & $-\infty$ \\
$-10$ & --- & $0.64597784$ & $-\infty$ & $-\infty$ \\
$-4$ & --- & $0.91468136$ & $-\infty$ & $-\infty$ \\
$-3$ & $1$ & $1$ & $-\infty$ & $-\infty$ \\
$-2$ & $\pi/\sqrt{8}$ & $1.11072073$ & $-\infty$ & $-\infty$ \\
$-5/3$ & $2/\sqrt{3}$ & $1.15470054$ & $-\infty$ & $-\infty$ \\
$-3/2$ & $3\pi/8$ & $1.17809725$ & $-\infty$ & $-\infty$ \\
$-4/3$ & $\pi\sqrt{75/512}$ & $1.20239047$ & $-\infty$ & $-\infty$ \\
$-1$ & $\sqrt{\pi/2}$ & $1.25331414$ & $-\infty$ & $-\infty$ \\
$-2/3$ & $\sqrt{128/75}$ & $1.30639453$ & $-\sqrt{6}$ & $-2.44948974$ \\
$-1/2$ & $4/3$ & $1.33333333$ & $-2$ & $-2$ \\
$-1/3$ & $\pi\sqrt{3}/4$ & $1.36034952$ & $-\sqrt{3}$ & $-1.73205081$ \\
$0$ & $\sqrt{2}$ & $1.41421356$ & $-\sqrt{2}$ & $-1.41421356$ \\
$1$ & $\pi/2$ & $1.57079633$ & $-1$ & $-1$ \\
$2$ & --- & $1.71731534$ & $-\sqrt{2/3}$ & $-0.81649658$ \\
$3$ & --- & $1.85407468$ & $-1/\sqrt{2}$ & $-0.70710678$ \\
$4$ & --- & $1.98232217$ & $-\sqrt{2/5}$ & $-0.63245553$ \\
$10$ & --- & $2.62843161$ & $-\sqrt{2/11}$ & $-0.42640143$ \\
$100$ & --- & $7.20340190$ & $-\sqrt{2/101}$ & $-0.14071951$ \\
$\infty$ & $\infty$ & $\infty$ & 0 & 0 \\
[0.5ex]
\bottomrule
\end{tabularx}
\caption{List of dimensionless collapse times $\tau$, given by \eqs{collapsetime}{collapsetimespecial}, and dimensionless collapse point velocities $\dot r(\tau)$, given by \eq{collapsedotr}. Exact solutions are given where they are algebraic, at least up to a factor $\pi$.}
\label{tab:tau}
\end{table}

\subsection{Noteworthy special solutions}

For completeness, I note that for select values of $\gamma$, \eq{radiussolution.a} can be reduced to well-known more compact solutions,
\begin{subnumcases}
	{r(t) = \label{eq:specialsolution}}
	\sqrt{1-t^2} & if $\gamma=-3$\label{eq:specialsolutionm3}\\
	1-t^2/2 & if $\gamma=0$\label{eq:specialsolution0}\\
	\cos(t) & if $\gamma=1$\label{eq:specialsolution1}\\
	{\rm cn}\big(t,\tfrac{1}{2}\big) & if $\gamma=3$\label{eq:specialsolution3},
\end{subnumcases}
where ${\rm cn}(x,y)$ is the Jacobi elliptic cosine function.

To my knowledge, Equations~(\ref{eq:specialsolution}a--c) are the only explicit closed-form solutions, \ie solutions in terms of commonly accepted basic functions. These three functions respectively describe an arc of a circle \eqp{specialsolutionm3}, a parabola \eqp{specialsolution0}, and a harmonic oscillation \eqp{specialsolution1}. Textbook examples of corresponding physical systems are given in \tab{applications}.

\eq{specialsolution3} is a special case of the Duffing equation \citep{Duffing1918}, describing an undamped and unforced anharmonic oscillator. General solutions in terms of the Jacobi elliptic family have recently been derived \citep{Salas2021}.

For some other values of $\gamma$, parametric solutions have been found well before explicit solutions. A famous example is the case of $\gamma=-2$, describing, \eg the spherical gravitational collapse (see \tab{applications}). In this case \citep{Lin1965},
\be\label{eq:parametric}
	t=\frac{\theta+\sin\theta}{\sqrt{8}},\quad r=\frac{1+\cos\theta}{2}.
\ee
Varying $\theta$ from $0$ to $\pi$ generates the collapse phase, whereas varying it from $-\pi$ to $0$ generates the symmetric growth phase. More lengthy parametric solutions have also been presented for other non-trivial cases, such as empty spherical cavitation bubbles ($\gamma=-4$) \cite{Mancas2016}.


\section{Discusion}\label{s:discussion}

The general solution of \eqs{main.standard}{main.ics.standard} given by \eq{radiussolution} warrants a brief discussion. As a preliminary remark, the linear transformation \eqp{nondimensional} between the normalised form $r(t)$ and its dimensional analogue $R(T)$ makes it straightforward to apply all properties of $r(t)$ to $R(t)$, and vice versa. For instance, as mentioned in \s{introduction}, any solution of \eqs{main}{main.ics} is invariant under time reversal. This symmetry equally applies to dimensionless coordinates, \ie $r(t)=r(-t)$, for all $t$, where a solution $r(t)$ exists. For efficiency, I will therefore limit the discussion in this section to the normalised form $r(t)$ and to positive times.

\subsection{Behaviour at collapse point}\label{ss:collapsepoint}

As illustrated in \fig{rt}, the steepness of the function $r(t)$ monotonically increases as the time sweeps from $t=0$ ($r=1$) to the collapse time $t=\tau$ ($r=0$). \eq{integral} shows this immediately and reveals that in the limit $t\rightarrow\tau_-$, the velocity $\dot r$ becomes
\be\label{eq:collapsedotr}
	\dot r(\tau)= 
	\left\{\begin{array}{ll}
		-\sqrt{2/(1+\gamma)} & \text{if }\gamma>-1\\
		-\infty & \text{if }\gamma\leq-1
	\end{array}\right..
\ee
Some explicit values of $\dot r(\tau)$ are listed in \tab{tau}.

The divergent behaviour of $\dot r(\tau)$ for $\gamma\leq-1$, which can often be interpreted as a positively diverging kinetic energy ($\propto\dot r^2$), is mimicked by a negatively diverging potential $\phi(r)$. In fact, following \eq{potential},
\be\label{eq:collapsephi}
	\lim_{r\rightarrow0_+}\phi(r)= 
	\left\{\begin{array}{ll}
		0 & \text{if }\gamma>-1\\
		-\infty & \text{if }\gamma\leq-1
	\end{array}\right.,
\ee
as is apparent in \fig{potential}.

\eqs{collapsedotr}{collapsephi} highlight the existence of two distinct regimes in the domain $\gamma\in\R$, separated by the critical value $\gamma=-1$. This value corresponds to the only singularity of \eq{radiussolution.a}, which is removed by \eq{radiussolution.b}. The dimensionless collapse time for $\gamma=-1$ is exactly $\sqrt{\pi/2}$ \eqp{collapsetimespecial}, the value shown by the dashed vertical line in \fig{rt}. Following \eq{collapsedotr}, all curves reaching $r=0$ to the left of this line do so vertically, whereas those to the right of the line come down at a finite slope.

The two regimes separated by $\gamma=-1$ correspond to different classes of physical problems (\cf \tab{applications}), explaining why---to my knowledge---they have never been addressed simultaneously in previous literature.

The regime $\gamma<-1$, sometimes including $\gamma=-1$, often describes a spherical collapse---gravitational or hydrodynamic in nature (see \tab{applications}). In this case, the collapse motion $r(t)$ near $t=\tau$, is sometimes referred to as `violent collapse' or `catastrophic collapse'. This description is quite literal, a saddening example being the implosion of the Titan submersible near the wreck of the Titanic in the North Atlantic Ocean (18~June 2023). This implosion was likely approximated by \eq{main} with $\gamma=-4$, predicting a collapse time $T_c\approx5\,\rm ms$, based on a capsule radius $R_0\approx1\,\rm m$, a driving pressure $\Delta p\approx350\,\rm bar$ and a water density $\rho\approx10^3\,\rm kg\,m^{-3}$.

The infinite velocity at the collapse point for $\gamma\leq-1$ is unphysical, as all real-world systems modelled by \eq{main} start deviating from this model near the collapse point. Secondary mechanisms, which might have been negligible for most of the collapse motion, suddenly become dominant, preventing the singularity. For example, in the case of collapsing cavitation bubbles, these mechanisms include shock waves (liquid compressibility), sonoluminescence, sonochemistry, vapour compression and micro-jetting \cite{Obreschkow2013c}. In the case of matter collapsing by self-gravity, the mechanisms preventing the divergence could be asymmetries, smoothing the collapse point, and pressure forces, possibly enhanced by strong radiation and phase transitions, \eg to a neutron star. If these mechanisms cannot prevent $\dot R(T)$ from approaching the speed of light, general relativistic effects take over, transforming the collapsing mass into a static black hole \citep{Oppenheimer1939}. This is likely the fate of the cores of massive stars ($\gtrsim25\,\msun$) \citep{Mirabel2017} and possibly also primordial gas clouds ($\sim10^5\msun$), collapsing directly into black holes \citep{Loeb1994}.

\subsection{Continuation past collapse point}\label{ss:domain}

The solution of \eq{radiussolution} is valid only on the interval $t\in[-\tau,\tau]$. However, \eq{main.standard} integrates smoothly past the collapse point $t=\tau$, to negative values of $r$, if $\gamma\in\N_0$.

For non-negative even $\gamma$, the sign of $\ddot r$ remains negative as $r<0$. Hence the `collapse' motion continues through the point $t=\tau$ with ever increasing velocity, such that $\lim_{t\rightarrow\infty}\dot r=-\infty$. The simplest example would be an object dropped to the ground (defined as $R=0$) into a vertical shaft ($R<0$) with constant acceleration $g$ ($\gamma=0$, see \tab{applications}).

For positive odd $\gamma$, however, $r(t)$ passes through $r(\tau)=0$, while inverting the sign of $\ddot r$. Symmetry considerations imply that $r(t)$ then describes an oscillating curve of wavelength $\lambda=4\tau$, made of reflections of the arc $r(t\in[0,\tau])$. This curve satisfies $r(0)=1$, $r(\tau)=0$, $r(2\tau)=-1$, $r(3\tau)=0$, $r(4\tau)=1$, and so forth. For $\gamma=1$, the oscillation is harmonic \eqp{specialsolution1}, for all larger odd $\gamma$, it is anharmonic (\eg Equation~\ref{eq:specialsolution3}).

Non-integer positive $\gamma$ generally lead to complex solutions $r(t)\in\C$, if integrated past the collapse point. In fact, immediately past this point, the acceleration has the complex argument $\arg(\ddot r)=(\gamma+1)\pi$, implying that its imaginary part only vanishes for integer $\gamma$.

For negative $\gamma$, \eq{radiussolution} diverges as $r\rightarrow0$. This is a coordinate singularity, which can, in principle, be removed through regularisation techniques, such 2-body and $N$-body regularisation in astrophysical gravitational simulations \citep{Aarseth2009}.

\begin{figure}
    \centering
    \includegraphics[width=\columnwidth]{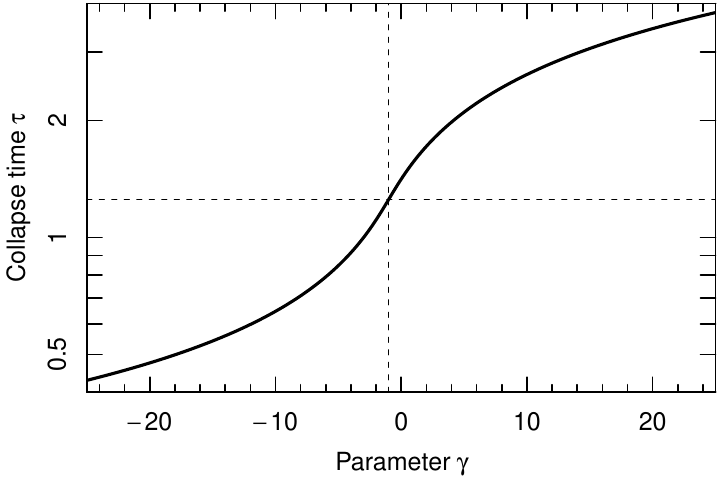}\vspace{-2mm}
    \caption{Dimensionless collapse time $\tau=T_c/T_0$ as a function of the control parameter $\gamma$. Plotted in semi-log coordinates, this function exhibits a non-trivial point symmetry about $(\gamma=-1,\tau=\sqrt{\pi/2})$, marked by the dashed lines. See \ss{hiddensymmetry}.}
    \label{fig:tau}
\end{figure}

\subsection{Collapse time symmetry}\label{ss:hiddensymmetry}

\fig{tau} shows the dimensionless collapse time $\tau$ as a function of the control parameter $\gamma$. With $\tau$ on a logarithmic axis, this function becomes a symmetric S-shape, which owes its symmetry to the beta function identity $B(\alpha,\beta)B(\alpha+\beta,1-\beta)=\pi/(\alpha\sin(\pi\beta))$. Applied to \eq{collapsetime}, this identity implies that any point $(\gamma,\tau(\gamma))$ can be mapped onto a corresponding point
\be\label{eq:gammasymmetry}
	(\gamma,\tau)\mapsto\left(\gamma'=-2-\gamma,\tau'=\frac{\pi}{2\tau}\right),
\ee
which also lies on the function $\tau(\gamma)$. This mapping is symmetric under the exchange of primed and non-primed variables. The only invariant point of this symmetry transformation, $\gamma=\gamma'=-1$ and $\tau=\tau'=\sqrt{\pi/2}$, is marked by the dashed lines in \fig{tau}. This point coincides with the special case of \eq{collapsetimespecial}, the singularity of \eq{radiussolution.a} and the transition of $\dot r(\tau)$ from a finite to a diverging value (\ss{collapsepoint}).


\begin{figure}
    \centering
    \includegraphics[width=\columnwidth]{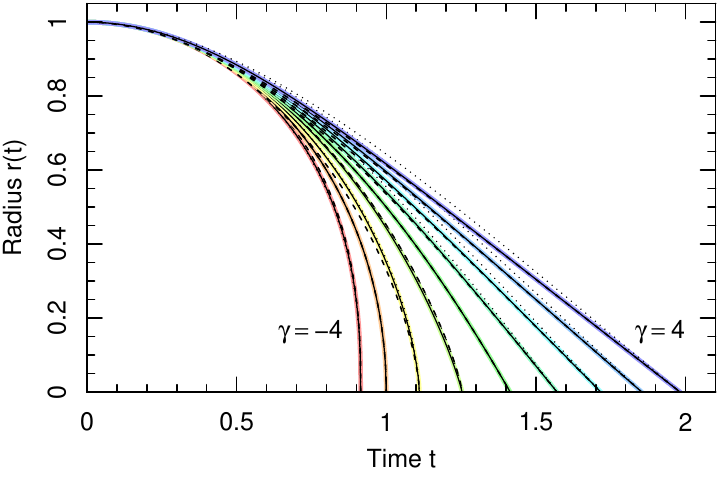}\vspace{-2mm}
    \caption{Exact versus approximate solutions of the normalised collapse equation \eq{main.standard} with initial conditions of \eq{main.ics.standard}, evaluated on the normalised time interval $t\in[0,\tau(\gamma)]$. The selected values of $\gamma$ are the same as in \fig{rt}. Solid lines are the exact solutions given in \eq{radiussolution}, whereas dashed and dotted lines use the approximation of \eq{approximation} with parameters $p_1$ and $q_1$ (dashed) and $p_2$ and $q_2$ (dotted).}
    \label{fig:approx}
\end{figure}

\subsection{General symmetry}\label{ss:generalsymmetry}

A noteworthy property of \eq{main.standard} is its formal invariance under a substitution $r\mapsto r':=r^\delta$, if $\delta=(1-\gamma)/2$, provable by invoking \eq{hybrid}. The restriction to collapse motions, implying $r,r'\leq1$ and $\ddot r,\ddot r'<0$, limits the applicability of this substitution to positive $\delta$, \ie to $\gamma<1$. A few elementary steps show that this substitution corresponds to the transformation
\be\label{eq:generalsymmetry}
	(\gamma,t,r)\mapsto\left(\gamma'\!=\frac{\gamma+3}{\gamma-1},t'\!=t\sqrt{\frac{1-\gamma}{2}},r'\!=r^{\frac{1-\gamma}{2}}\right)\!.
\ee
Like \eq{gammasymmetry}, this mapping is symmetric under the exchange of primed and non-primed variables; and likewise this symmetry transformation exhibits a single invariant point at $\gamma=\gamma'=-1$, where $(r',t')=(r,t)$.

The symmetry of \eq{gammasymmetry} can help identify new solutions based on existing ones, an example being the derivation of \eq{specialsolutionm3} from \eq{specialsolution0}, and vice versa.

When applying \eq{generalsymmetry} to dimensional quantities $R$ and $T$ instead of $r$ and $t$, the physical constant $k$ \eqp{main} must also be transformed, as $k'=k\delta R_0^{1+\gamma}$. For example, Rayleigh's \eq{rayleigh} can be written as
\be\label{eq:rayleigh2}
	\frac{\d^2Z}{\d t^2} = -\frac{5\Delta p}{2\rho}Z^{1/5},
\ee
where $Z:=R^{5/2}$.

\subsection{Polynomial approximations}

In some practical cases, it may be useful and sufficient to substitute the exact solution $r(t)$ given in \eq{radiussolution} by simple approximations $\tilde r(t)$. A possible choice is
\begin{subnumcases}
	{\hspace{-5mm}\tilde r(t) = \label{eq:approximation}}
	\!\left(1-x^2\right)^p & \!\!\!if $\gamma\leq-1$ \label{eq:approximation.a}\\
	\!q(1\!-\!|x|)\!-\!(q\!-\!1)(1\!-\!|x|)^{\frac{q}{q-1}} & \!\!\!if $\gamma>-1$\label{eq:approximation.b},
\end{subnumcases}
where $x=t/\tau$ (with $\tau$ given in \eqs{collapsetime}{collapsetimespecial}) and where $p\in(0,1)$ and $q>1$ are shape parameters.

These are arguably the simplest polynomials, which simultaneously satisfy six essential properties of the exact solution: (1) time symmetry, $\tilde r(t)=\tilde r(-t)$ if $t\in[-\tau,\tau]$; (2) negative curvature $\ddot{\tilde r}(t)<0$ if $t\in(-\tau,\tau)$, (3) $\tilde r(0)=1$; (4) $\tilde r(\tau)=0$; (5) $\dot{\tilde r}(0)=0$; and (6) $\lim_{t\rightarrow\tau}\dot{\tilde r}(t)$ is finite if and only if $\gamma>-1$.

Let us first consider \eq{approximation.a} (case $\gamma\leq-1$), where the exponent $p$ needs to be specified. Choosing $p$ equal to $p_1=2/(1-\gamma)$ (only for $\gamma$ strictly below $-1$) ensures the correct asymptotic behaviour of $\dot{\tilde r}$ as a function of $\tilde r$ in the limit $t\rightarrow\tau$. Earlier, I introduced this approximation specifically for $\gamma=-4$ (thus $p=2/5$) \cite{Obreschkow2012a}. \eq{approximation.a} can then be seen as the first term in a rapidly converging series, which tends toward the true solution if using exact analytical coefficients \cite{Amore2013}.

Alternatively, setting $p$ equal to $p_2=\tau^2/2$ (for all $\gamma\leq-1$), satisfies $\ddot{\tilde r}(0)=\ddot r(0)=-1$. For $\gamma=-2$, the case of the gravitational collapse, this choice is $p=(\pi/4)^2\approx0.6169$. This is nearly identical to an analogous recent approximation \cite{Girichidis2014}, which reads $r(t)=(1-t^2)^p$ with $p=1.8614/3\approx0.6205$ determined via a fitting technique.

For \eq{approximation.b} (case $\gamma>-1$), we have similar options for the parameter $q$. Setting it to $q_1=\eta B(\eta,\tfrac{1}{2})$ (with $\eta$ defined in Equation~\ref{eq:defeta}) ensures that $\dot{\tilde r}(\tau)=\dot r(\tau)$. Alternatively, setting it equal to $q_2=\tau^2/(\tau^2-1)$ ensures $\ddot{\tilde r}(0)=\ddot r(0)=-1$.

As shown in \fig{approx}, both of the above choices for $p$ and $q$ provide reasonable approximations of the exact solution. In general, the best pick depends on the application. If $\gamma=-3$ (hence $p_1=p_2=\tfrac{1}{2}$) or $\gamma=0$ (hence $q_1=q_2=2$), the two options are identical and \eqs{approximation.a}{approximation.b} become the exact solutions of \eqs{specialsolutionm3}{specialsolution0}, respectively. Interestingly, these two solutions are related to each other via the symmetry transformation of \eq{generalsymmetry}.

A more quantitative discussion of the approximations is beyond the scope of this paper, but can be found elsewhere for $\gamma=-4$ \citep{Obreschkow2012a,Amore2013}.


\section{Conclusion}\label{s:conclusion}

This paper investigated the differential \eq{main} with boundary conditions stated in \eq{main.ics}. I have called this the spherical collapse equation in reference to its most common applications (\cf \tab{applications}).

I have shown that this equation exhibits a unified explicit solution, invoking the beta distribution, a fundamental probability density function in statistics. Elements of the derivation can be found elsewhere, spread over decades of literature across several unrelated fields. The contribution of this work is to unify previously known parametric, implicit and explicit solutions from hydrodynamics (\cf Rayleigh-Plesset equation), hydrostatics (Lane-Emden equation) and astrophysics (top-hat spherical collapse model), all limited to different subsets of $\gamma$, into a single concise solution (Equation~\ref{eq:radiussolution.a}), universally valid for all real $\gamma\ne-1$. For $\gamma=-1$ this solution is singular, but its limit $\gamma\rightarrow-1$ exists, allowing the singularity to be removed via \eq{radiussolution.b}.

This general character ($\gamma\in\R$) of the explicit solution unifies seemingly unrelated physical applications. In this regard, this work provides context for recently found apparent analogies between cavitation bubbles and processes in astronomy and cosmology \citep{Farhat2020,Rousseaux2020,Rosu2021,Shneider2021}.

The existence of a general explicit solution of \eq{main} is barely known, particularly in astrophysics and cosmology, where only parametric solutions like \eq{parametric} are commonly taught. 
Not only is the solution of \eq{radiussolution} more elegant and faster to evaluate than numerical integration, but its explicit form can serve as a basis for further analytical manipulations and derivations to analyse particular problems. Illustrating this point, the theory of galaxies has benefited from Freeman's analytic solution for the circular velocity of exponential disks \citep{Freeman1970}. This velocity can also be computed by direct numerical integration, but the explicit expression in terms of Bessel functions has greatly simplified further studies.

The explicit solution also offers insights pertaining to symmetry properties (Sections~\ref{ss:hiddensymmetry} and \ref{ss:generalsymmetry}) and a connection between spherical collapse motions and probability theory. Finally, \eq{radiussolution} can serve as an exact benchmark for testing numerical integration techniques and simulation codes, \eg for hydrodynamic and gravitational simulations.

An avenue for research in pure mathematics is an extension to complex-valued functions $r(t)$, where additional symmetries in the complex plane exist for $\gamma\in\N$.


\begin{acknowledgments}
I thank Dr.~Martin Bruderer, Dr.~Aaron Ludlow, Dr.~Zachari Slepian, Dr.~Dmitry Sinelshchikov, Dr.~Roberto Iacono and Prof.~Stefani Mancas for insightful suggestions that have inspired and improved this work. I also acknowledge both anonymous referees for their comments and positivity, as well as an anonymous referee of an earlier draft, for pointing out the substitutions discussed in \ss{generalsymmetry}. Last but not least, I extend my gratitude to Dr.~Mohamed Farhat and Prof.~Andrea Prosperetti for many past discussions, which were a source of inspiration in writing this paper. I am a recipient of an Australian Research Council Future Fellowship (FT190100083) funded by the Australian Government.
\end{acknowledgments}


\appendix

\section{Analytical derivation}\label{a:derivation}

\eq{integral.a} can be rewritten to separate time from position coordinates,
\be
	\d t^2 = \frac{1+\gamma}{2(1-r^{1+\gamma})}\d r^2,
\ee
During the collapse ($0<r<1$), the RHS is positive for any $\gamma\ne-1$. Since the collapse is characterised by a shrinking radius ($\d r<0$) with growing time ($\d t>0$), we are interested in the negative branch,
\be
	\d t = -\left[\frac{1+\gamma}{2(1-r^{1+\gamma})}\right]^{1/2}\d r.
\ee
Let us integrate this equation backward in time from the dimensionless collapse time $\tau$ to an earlier time $t>0$,
\be
	\int_\tau^t \d t = -\left[\frac{|1+\gamma|}{2}\right]^{1/2}\int_0^{r(t)}\frac{\d r}{|1-r^{1+\gamma}|^{1/2}}.
\ee
The lower bound of the second integral is zero by definition of the collapse time. To simplify the notation of the following equations let us introduce the positive constant
\be
	\eta=\frac{1}{|1+\gamma|},
\ee
and let $s:=r^{|1+\gamma|}\in[0,1]$. With these substitutions, $r=s^\eta$ and $\d r=\eta s^{\eta-1}\d s$, and hence
\be
	\int_\tau^t \d t = -\sqrt{\frac{\eta}{2}}~\int_0^{s(t)}\frac{s^{\eta-1}\d s}{|1-s^{\pm1}|^{1/2}},
\ee
where the sign of the exponent in the denominator is equal to the sign of $1+\gamma$. If this sign is negative, we multiply the numerator and denominator of the fraction by $s^{1/2}$. Then,
\be\label{eq:timeintegral}
	\int_\tau^t \d t = -\sqrt{\frac{\eta}{2}}~\int_0^{s(t)}\frac{s^{\alpha-1}\d s}{(1-s)^{1/2}},
\ee
with $\alpha=\eta$, if $\gamma>-1$, and $\alpha=\eta+\tfrac{1}{2}$, if $\gamma<-1$. We can readily unify these two cases in the single equation
\be
	\alpha = \frac{1}{4}+\frac{3-\gamma}{4|1+\gamma|}.
\ee
We have defined $\alpha$ in this way to make the right-hand integral of \eq{timeintegral} identical to the definition of the \emph{incomplete beta function}, $B(x;\alpha,\beta)=\int_0^xy^{\alpha-1}(1-y)^{\beta-1}\d y$. Hence, \eq{timeintegral} solves to \eq{timesolution1}.


\section{Code examples}\label{a:implementation}

The quantile function $Q(x;\alpha,\beta)$, which forms the heart of the explicit \eq{radiussolution.a}, is readily accessible in most programming languages: \code{qbeta} in \Rcode, \code{scipy.stats.beta.ppf} in \Pythoncode, \code{BETA.INV} in \textsc{Excel}, \code{InverseBetaRegularized} in \textsc{Mathematica}, \code{InvIncompleteBeta} in the \textsc{ALGLIB} library for \textsc{C++}, \textsc{C\#}, \textsc{Java}, \textsc{Python}, \textsc{Delphi}; etc. For reference, this section provides a few explicit code examples for plotting the collapse functions shown in \fig{rt}.

\subsection{Gravitational collapse}

The gravitational collapse of a uniform pressure-free sphere is governed by \eq{tophat}. The solution $r(t)$ in natural units $R_0$ and $T_0=\sqrt{R_0^3/(GM)}$ is given by \eq{radiussolution.a} with $\gamma=-2$. Perhaps the most compact implementation is achieved in the statistical language \Rcode, where three lines suffice to evaluate and plot this explicit solution:

\definecolor{bg}{rgb}{0.9,0.9,0.9} 
\lstset{
backgroundcolor=\color{bg},  
showstringspaces=false  
}

\begin{lstlisting}[language=R,basicstyle=\small]
tau = beta(3/2,1/2)/sqrt(2)
r = function(x) qbeta(1-x/tau,3/2,1/2)
curve(r,0,tau,1000)
\end{lstlisting}
The resulting plot is analogous to the yellow line for $\gamma=-2$ in \fig{rt}. Most astrophysicists might be more accustomed to \Pythoncode, where this code could read:

\begin{lstlisting}[language=Python,basicstyle=\small]
import numpy as np
import matplotlib.pyplot as plt
from scipy.stats import beta
from scipy.special import beta as betafct

tau = betafct(3/2, 1/2)/np.sqrt(2)
def r(t):
  return beta.ppf(1-t/tau, 3/2, 1/2)

tvalues = np.linspace(0, tau, 1000)
plt.plot(tvalues, r(tvalues))
plt.show()
\end{lstlisting}

\subsection{Cavitation bubble collapse}

Readers closer to cavitation dynamics than astrophysics may consider the example of a spherical cavity collapsing in an incompressible liquid of density $\rho$, without viscosity and surface tension. Assuming a constant driving pressure $\Delta p$, the bubble evolution is governed by \eq{rayleigh}; hence the collapse motion $r(t)$ in natural units $R_0$ and $T_0=R_0\sqrt{\rho/\Delta p}$ is given by \eq{radiussolution.a} with $\gamma=-4$. The corresponding \Rcode script reads:

\begin{lstlisting}[language=R,basicstyle=\small]
tau = beta(5/6,1/2)/sqrt(6)
r = function(x) qbeta(1-x/tau,5/6,1/2)^(1/3)
curve(r,0,tau,1000)
\end{lstlisting}
This code reproduces the red line for $\gamma=-4$ in \fig{rt}. In \Pythoncode, this could be implemented as:

\begin{lstlisting}[language=Python,basicstyle=\small]
import numpy as np
import matplotlib.pyplot as plt
from scipy.stats import beta
from scipy.special import beta as betafct

tau = betafct(5/6, 1/2)/np.sqrt(6)
def r(t):
  return beta.ppf(1-t/tau, 5/6, 1/2)**(1/3)

tvalues = np.linspace(0, tau, 1000)
plt.plot(tvalues, r(tvalues))
plt.show()
\end{lstlisting}


%

\end{document}